\documentstyle[11pt]{article}
\newcommand{\calp}{{\cal P}}            
\newcommand{\calb}{{\cal B}}            
\newcommand{\calep}{{\cal EP}}            

\newcommand{\ind}{{\bf 1}}

\newcommand{\startb}{\parindent0pt\bf}  

\newtheorem{thm}{Theorem}
\newtheorem{lem}[thm]{Lemma}

\newtheorem{defn}{Definition}
\newtheorem{rem}{Remark}

\begin{document}

\baselineskip18pt

\title{{\Large\bf Linear Optimization over a Polymatroid with Side
Constraints \\
-- Scheduling Queues and Minimizing Submodular Functions }}
\author{Yingdong Lu\\
IBM T.J.\ Watson Research Center\\ Yorktown Heights, N.Y.10598 \\ \\
David D.\ Yao\\ IEOR Dept, 302 Mudd Building
\\Columbia University, New York, NY 10027}
\date{April 2007}
\maketitle

\begin{abstract}
Two seemingly unrelated problems, scheduling a multiclass
queueing system and minimizing a submodular function,
share a rather deep connection via the polymatroid that is characterized
by a submodular set function on the one hand and represents
the performance polytope of the queueing system on the other
hand. We first develop what we call a {\it grouping} algorithm
that solves the queueing scheduling problem under side
constraints, with a computational effort of $O(n^3LP(n))$,
$n$ being the number of job classes, and $LP(n)$ being the computational 
efforts of solving a linear program with no more than $n$ variables and $n$ 
constraints.
The algorithm organizes the job classes into
groups, and identifies the optimal policy to be a priority rule
across the groups and a randomized rule within each group
(to enforce the side constraints).
We then apply the grouping algorithm
to the submodular function minimization,
mapping the latter to a queueing scheduling problem with
side constraints.
We show the minimizing subset can be identified by
applying the grouping algorithm $n$ times.
Hence, this results in a algorithm
that minimizes a submodular function with an effort of $O(n^4LP(n))$.
\end{abstract}

\section{Introduction}
\label{sec:intro}


Linear optimization over a polymatroid --- a polytope that
possesses matroid properties --- is known to be solvable through
a greedy algorithm; refer to
Edmonds \cite{edmonds}, Dunstan and Welsh
\cite{dunstan},  and Welsh \cite{welsh}.
It is also known that the performance space of a
multiclass queueing system is a polymatroid, provided
the system satisfies strong conservation laws; refer to
Shanthikumar and Yao \cite{syor}, and also the earlier works
of Coffman and Mitrani \cite{coffman}, and Federgruen and Groenevelt
\cite{fgor,fgmgc}.
With this connection,
the control or dynamic scheduling of multiclass queues
can be translated into an optimization problem; and under a
linear objective, the optimal control is readily identified
as a priority rule associated with one of the vertices of
the performance polytope.

In many scheduling applications, there are additional,
``side constraints'' on the performance; for instance, the
delay of jobs should not exceed a certain limit, or the
completion rate (throughput) of jobs must meet a certain requirement.
The dynamic scheduling of multiclass queues under
such side constraints is still polynomially solvable,
as argued in Yao and Zhang
\cite{yaozhang}, directly applying a result by Gr\"otschel,
Lov\'asz and Shrijver \cite{GLS81}.
This argument, however, is essentially an ``existence proof.''
It does not explicitly identify an optimal policy, let alone
one with a simple structure that appeals to implementation.


In this paper, we develop a combinatorial algorithm that identifies
the optimal policy for scheduling multiclass queues with side
constraints. This is what we call a ``grouping algorithm,'' which
partitions the $n$ job classes into groups, and the associated
optimal policy is a priority rule across the groups, while enforcing
the side constraints (via randomization, for instance) among the
classes within each group. The grouping algorithm identifies the
optimal policy with a computational effort of $O(n^3\times LP(n))$,
where $LP(n)$ refers to the complexity of solving a linear program
with number of variables and constraints bounded by $n$.

It turns out that this grouping algorithm can also be used
to solve the problem of minimizing a submodular function
--- finding the subset that minimizes a set function,
$\psi:2^E\mapsto \cal R$, which maps the subsets of a given set
$E$ to the real line and satisfies submodularity
(refer to Definition \ref{def:prhs} (iii)).
Here the key ingredient is that
a certain representation of the submodular function minimization
can be turned into a queueing scheduling problem with side constraints,
with the submodular function relating to the polymatroid that
is the performance polytope of the queueing system.
We show that the minimizing subset can be identified by
applying the grouping algorithm
$n=|E|$ times; each time, it identifies
a subset, and only these subsets
are candidates for optimality.
Hence, the overall computational effort is $O(n^4LP(n))$.

Submodular function minimization (SFM), widely regarded as  ``the most
important problem related to submodular function''
(\cite{glsbook}), has attracted much research
effort in recent years.
Here we briefly highlight the relevant literature, as
the details are now available in two recent survey papers by
Fleischer \cite{fleischer2} and McCormick \cite{mccormick}.

Using the ellipsoid approach,
Gr\"otschel, Lov\'asz and Shrijver
\cite{GLS81} established that the SFM is polynomially solvable.
However, the ellipsoid technique is not usually
a practical algorithm, and it does not generate much
combinatorial insight.
Cunningham \cite{cunningham} proposed a max-flow type of
approach based on linear programming (LP), which results in a
combinatorial algorithm that solves the SFM in pseudo-polynomial time.
(The qualifier ``pseudo'' refers to the fact that the
running time is not only a polynomial in $n$, but also a
polynomial in $M$, an upper bound of the set function.)
Recently, combinatorial algorithms that solves SFM in
strongly polynomial time (meaning, the running time is
a polynomial in $n$ only, and is independent of $M$)
have been independently developed  by Iwata, Fleicher and
Fujishige \cite{IFF}, and by Shrijver \cite{Shrijver}; and
further improvement has appeared in Iwata \cite{iwata,iwata1},
and Fleischer and Iwata \cite{fi}.
The starting point of these algorithms is
an equivalent  problem over
a polymatroid with side constraints, and this LP is solved
via an augmenting-path approach.
In contrast, our approach to the SFM is to solve a sequence of $n$ LP
problems;  each is an LP over a polymatroid with side constraints,
but the inequalities in the side constraints are in the opposite direction
to those in existing SFM algorithms. None of these $n$ LP's is equivalent
to the SFM, but solving all of them does yield the solution to the
SFM. More importantly, each such LP can be represented as a
queueing scheduling problem highlighted above,
and hence can be solved by the grouping algorithm.
(Refer to more details in \S\ref{sec:minsubmod}, Remarks \ref{rem:ours}
and \ref{rem:others} in particular.)
The advantage of our algorithm is that it is fully combinatorial
(meaning, it only uses addition, subtraction and comparison);
and its running time of $O(n^4LP(n))$ matches that of the strongly
polynomial version of the ellipsoid algorithm, and is significantly better
than the fastest existing combinatorial algorithms.
(For instance, the one due to Iwata \cite{iwata1} has a running time of $O(n^7)$;
refer to \cite{mccormick}.)


Briefly, the rest of the paper is organized as follows.
In the next section, we review preliminary materials
for linear optimization over a polymatroid, and its
connection to ququeing scheduling via strong conservation laws.
In \S\ref{sec:algorithm} we collect the essentials of the grouping algorithm that
solves the queueing scheduling problem with side constraints, in particular
the conditions that will ensure its feasibility and optimality.
This is folllowed by a detailed description of the algorithm in
\S\ref{sec:algorithm1}, where we also show that the required
conditions outlined in \S\ref{sec:algorithm} are all met.
In \S\ref{sec:minsubmod}
we apply the grouping algorithm to the SMF problem.
Finally, in \S\ref{sec:symmetry}
we consider a special case when the submodular function satisfies the
so-called generalized symmetry.

The polymatroid connection to queueing scheduling
has in recent years been extended to a more general framework
via so-called
generalized conservation laws;
refer to the papers by Bertsimas and
his co-workers, e.g., \cite{bertsimquesta,BNM,mitside}; also
refer to Lu \cite{lyd}, Zhang \cite{li},
and Chapter 11 of Chen and Yao \cite{cynet}.
In this paper,  however, we choose to
stay within the more restrictive framework of polymatroid,
since our emphasis here is not on queueing scheduling per se,
but rather on the interplay between such scheduling problems under
side constraints and minimizing submodular functions.

\section{Preliminaries}
\label{sec:prelim}

\subsection{LP over a Polymatroid}
\label{sec:lppoly}

Throughout, $E=\{1,\dots, n\}$ is a finite set;
$x=(x_i)_{i\in E}$ is a  vector, and $x(A):=\sum_{i\in A} x_i$
for any $A\subseteq E$.
We shall use ``increasing'' and ``decreasing'' in the non-strict
sense, to mean non-decreasing and non-increasing, respectively.

\begin{defn}
{\rm
(Welsh \cite{welsh}) Given $E=\{1,2,...,n\}$, the following polytope,
\label{def:prhs}
\begin{eqnarray}
\label{defPf}
\calp(f)=\{ x\ge 0: x (A) \le f(A), A \subseteq E\},
\end{eqnarray}
is called a polymatroid if the function $f:2^E \rightarrow {\cal
R}_+$ satisfies the following properties:

\noindent (i)  ({\it normalized}) $f(\emptyset) =0$;

\noindent (ii)  ({\it increasing}) if $A \subseteq B \subseteq E$,
then $f(A) \le f(B)$;

\noindent (iii)  ({\it submodular}) if $A, B \subseteq E$, then
$f(A)+ f(B) \ge f(A\cup B) +f(A \cap B)$
}
\end{defn}

Let $\pi$ be a permutation of $\{1,2,...,n\}$, and define a vector
$x^{\pi}$ with the following components:
\begin{eqnarray}
x_{\pi_1}^{\pi} &=&f(\{\pi_1 \}), \nonumber \\ x_{\pi_2}^{\pi}
&=& f(\{\pi_1,\pi_2\}) - f(\{\pi_1 \}),\nonumber \\
&\cdots&\label{vertex}\\ &\cdots&\nonumber\\ x_{\pi_n}^{\pi}
&=&f(\{\pi_1, \pi_2,...,\pi_n\}) - f(\{ \pi_1,
\pi_2,...,\pi_{n-1}\}).\nonumber
\end{eqnarray}
Then, we have another definition for polymatroid:

\begin{defn} 
\label{def:pvertex}
{\rm
The polytope $\calp(f)$ in (\ref{defPf}) is a polymatroid
if $x^{\pi} \in \calp(f)$ for any permutation $\pi$.
}
\end{defn}

That the two definitions, \ref{def:prhs} and \ref{def:pvertex},
are equivalent is a classical result; refer to, e.g., Edmonds
\cite{edmonds}, Welsh \cite{welsh}, and Dunstan and Welsh
\cite{dunstan}. The following theorem, also a classical result,
can be found in these references as well.


\begin{thm}
\label{thm:plp}
{\rm
(\cite{dunstan,edmonds,welsh}) Consider the
LP problem: $\max_{x\in\calp (f)}\sum_{i\in E}c_ix_i$. Suppose
$\pi=(\pi_i)_{i\in E}$ is a permutation, such that $$c_{\pi_1} \ge
c_{\pi_2}\ge \cdots \ge c_{\pi_n}.$$ Then, $x^\pi$ as specified in
(\ref{vertex}) is the optimal solution.
}
\end{thm}

For the polymatroid $\calp(f)$, we call the following polytope its
{\it base}:
\begin{eqnarray}
\label{basef}
\calb(f)=\{ x\ge 0: \, x(A) \le f(A),\, A \subset E;
\, x(E)=f(E)\}.
\end{eqnarray}
%

In queueing applications, it is useful to have another
type of base, in addition to the one in (\ref{basef}):
\begin{eqnarray}
\label{defbase}
{\cal B}(b):=\{x\ge 0:\, x(A) \ge b(A),\, A
\subset E; \, x(E) = b(E)\},
\end{eqnarray}
where $b(\cdot)$ is increasing and supermodular, with
$b(\emptyset)=0$.
It is straightforward to verify that $\calb (b)=\calb (f)$ by
letting $f(A)=b(E)-b(E-A)$. Since the resulting $f$ is submodular,
increasing, and $f(\emptyset)=0$,
we know ${\cal B}(b)$ is also the base of a polymatroid.
In particular, each of its
vertices can be expressed as $x^\pi$ as in (\ref{vertex}), with
$f(\cdot)$ replaced by $b(\cdot)$.


\subsection{Polymatroid Structure in Multiclass Queues}
\label{sec:qnet}

Let $E=\{1,2,...,n\}$ denote the set of all job
classes (or, types) in a queueing system.
Let $x=(x_i)_{i\in E}$ denote the performance vector, e.g.,
waiting time (delay) or  service completions
(throughput). Let ${\cal U}$ denote the family of all {\it
admissible} policies that control the way jobs are served in
the system. This includes all policies that are non-anticipative (i.e., no
clairvoyance) and non-idling (i.e., serve whenever there is a job
to be served). We have:

\begin{defn}
\label{def:cl}
{\rm
(Shanthikumar and Yao \cite{syor}) $\;$ A
performance vector $x$ of a stochastic system with a set of
multiple job classes, $E=\{1,...,n\}$, is said to satisfy
conservation laws, if there exists a set function $b$ (resp.\
$f$): $2^E \mapsto {\cal R}_+$, satisfying
\begin{eqnarray}
\label{cfb} b(A) = \sum_{i\in A}x_{\pi_i}, \quad\forall
\pi:\{\pi_1,...,\pi_{|A|}\} = A, \quad\forall A \subseteq E;
\end{eqnarray}
or respectively,
\begin{eqnarray}
\label{cff} f(A) = \sum_{i\in A}x_{\pi_i}, \quad\forall
\pi:\{\pi_1,...,\pi_{|A|}\} =A, \quad\forall A\subseteq E;
\end{eqnarray}
such that under any admissible control,  $u \in {\cal U}$,
 the following is satisfied:
\begin{eqnarray*}
\sum_{i\in A}x^u_{\pi_i} \ge b(A), \quad\forall A \subset E;
\qquad\sum_{i \in E}x^u_{\pi_i} =b(E);
\end{eqnarray*}
or respectively,
\begin{eqnarray*}
\sum_{i\in A}x^u_{\pi_i} \le f(A), \quad\forall A \subset E; \qquad
\sum_{i\in E}x^u_{\pi_i} =f(E).
\end{eqnarray*}
}
\end{defn}

Intuitively, the above definition says:
\begin{itemize}
\item
The total performance over all job classes in $E$ is invariant: it
is a constant that is policy independent.

\item
The total performance over any given subset, $A \subset E$, of job
types is minimized (e.g., if $x$ is delay) or maximized (e.g., if
$x$ is throughput) by offering priority to job types in $A$ over
all other types in $E-A$.

\item
The function $b(A)$ or $f(A)$ is the corresponding performance
when $A$ is given priority over $E-A$. Furthermore, the
permutation $\pi$ represents a priority rule such that $\pi_1$ is
given the highest priority; $\pi_2$, the second highest priority;
and so forth. And the requirement, $\forall\pi:
\{\pi_1,...,\pi_{|A|}\} = A$, signifies that $b(A)$ and $f(A)$
depend only on the composition of the job types in $A$; they are
independent of the priorities among these job types.
\end{itemize}

The above definition of conservation laws, along with Definition
\ref{def:pvertex}, immediately leads to the important fact that
when conservation laws are satisfied, the performance space over
all admissible controls is a polymatroid --- in fact, the base
polytope of a polymatroid, and the
right-hand-side function is determined by
the conservation laws too, i.e. the performance resulting from
assigning priority to certain subset of jobs. This is simply
because any $\pi$, a priority rule, is admissible, and hence,
$x^\pi$ belongs to the performance polytope, and $x^\pi$ following
(\ref{cfb}) and (\ref{cff}) are of the same form as $x^\pi$ in
(\ref{vertex}).

Suppose we want to find the optimal control --- a dynamic
scheduling rule --- of the $n$ job types in the set $E$, so as to
maximize an objective function with respect to a performance vector $x$.
The problem can be expressed as follows:
\begin{eqnarray}
\label{qsched}
\max_{u\in {\cal U}} \sum_{i\in E}c_ix^u_i,
\end{eqnarray}
with $c_i$'s being the costs. Suppose $x$  satisfies conservation
laws; and hence, the performance space under all admissible
controls is the base polytope of a polymatroid, denoted $\calb
(f)$. So we can first solve a linear programming problem:
\begin{eqnarray*}
\max_{x \in \calb (f)} \sum_{i\in E}c_ix_i.
\end{eqnarray*}
From Theorem \ref{thm:plp}, we know the optimal solution is
$x^{\pi^*}$, where $\pi^*$ is a permutation that orders the cost
coefficients in decreasing order:
$$c_{\pi^*_1}\ge c_{\pi^*_2}\ge
\cdots \ge c_{\pi^*_n}.$$
We can then identify a control $u^*$
that realizes this performance. The fact that the solution being a
vertex of the performance space ensures that the optimal policy is
a priority rule.
It is clear
that to derive the optimal policy, there is no need to compute the
right-hand-side function of the performance polytope.
For more details regarding conservation laws and the connection
to polymatroid, including many examples, refer to Chapter 11 of
Chen and Yao \cite{cynet}.

\section{Side Constraints: Feasibility and Optimality}
\label{sec:algorithm}

%
Consider the polymatriod $\calp (f)$ in (\ref{defPf}).
Suppose there are additional constraints: $x_i \ge d_i\ge 0$, $i\in E$.
(For instance, $d_i$ is the requirement on the throughput of type $i$ jobs.)
Adding these to $\calp(f)$ results in the
following polytope:
\begin{eqnarray}
\label{polyside}
\calp'(f) = \{ 
x(A) \le f(A), A\subseteq E;\,
 x_i \ge d_i , i\in E\}
\end{eqnarray}
We want to solve the following LP problem:
\begin{eqnarray}
\label{lpside}
\max_{x\in\calp'(f)}  \sum_{k=1}^n c_k x_k.
\end{eqnarray}

For the LP in (\ref{lpside}) to be feasible,
the following condition is readily verified.

\begin{lem}
\label{lemrossyao}
{\rm
(Ross and Yao \cite{rossyao1}) The polytope
$\calp'(f)$ is non-empty if and only if $d:=(d_1,...,d_n )\in
\calp(f)$.
}
\end{lem}

\medskip

Below, we shall assume that $d\in \calp (f)$ always holds.
A direct application of the
following result establishes that
the LP in (\ref{lpside}) is polynomially solvable.

\begin{thm}
\label{thm:gls}
{\rm
(Gr\"otschel, Lov\'asz, and Schrijver \cite{GLS81})
$\;$ Let ${\cal K}$ be a class of convex bodies.
There exists a polynomial algorithm to solve
the separation problem for members of ${\cal K}$,
if and only if there exists a polynomial
algorithm to solve the optimization problem for the members of ${\cal K}$.
}
\end{thm}

Making use of this result, we know that the membership test problem
(which is equivalent to the separation problem) over $\calp (f)$ is polynomially
solvable, since the LP over $\calp (f)$ is.
Hence, the membership test over $\calp' (f)$ is also polynomially
solvable, since this amounts to checking the additional $n$
side constraints.
Therefore, the LP over $\calp' (f)$ is also polynomially solvable.

However, the argument in \cite{GLS81} is based on an ellipsoid approach.
In contrast, below we develop a fully combinatorial algorithm that solves the
LP over $\calp' (f)$ in $O(n^4)$ steps, and also reveals the structure of
the optimal solution.


\medskip

Without loss of generality, assume
 the coefficients in (\ref{lpside}) are ordered as follows:
\begin{eqnarray}
\label{corder}
c_1\ge c_2\ge \cdots \ge c_n >0.
\end{eqnarray}

The algorithm below divides the $n$ variables into $m\le n$ groups:
${G_1,...,G_m}$.
Denote, $G^{(0)}=\emptyset$; and
$$G^{(i)} := G_1\cup G_2\cup \cdots \cup G_i, \qquad i =1,...,m.$$
Within each group $G_i$, there is exactly one ``lead'' variable,
indexed $\ell_i$, which takes
on a value $x_{\ell_i}\ge d_{\ell_i}$;
the rest of the group, denoted $G_i^S$ (possibly
empty), consists of variables that all take values {\it at} the constraints,
i.e.,
\begin{eqnarray}
\label{nonleadvar}
x_{k_i}=d_{k_i}, \quad k_i\in G_i^S, \quad i=1,...,m.
\end{eqnarray}
Furthermore, the lead variable takes the following value:
\begin{eqnarray}
\label{leadvar}
x_{\ell_i}&=& f(G^{(i)})-f(G^{(i-1)})-d(G^S_i), \quad i=1,...,m ;
\end{eqnarray}
where we have used the notation $d(A):=\sum_{i\in A} d_i$, for
$A=G^S_i$.

\begin{lem}
\label{lem:feasible}
{\rm
The solution $x$ as specified in (\ref{nonleadvar}) and (\ref{leadvar})
is feasible, i.e., $x\in \calp' (f)$, if for each $i=1,...,m$,
we have: (a) $x_{\ell_i}\ge d_{\ell_i}$, and (b)
\begin{equation}
\label{g1a}
x(A) \le  f(G^{(i-1)};A) -f(G^{(i-1)}), \quad A\subset G_i .
\end{equation}
When the above do hold, we have
\begin{eqnarray}
x(G^{(i-1)};A) \le f(G^{(i-1)};A),
\quad A\subset G_i ;
\qquad
x(G^{(i)}) &=& f(G^{(i)}) ;
\label{g2}
\end{eqnarray}
for all $i=1,...,m$.
}
\end{lem}

{\startb Proof.}
First, all side constraints are satisfied, given (a), along with
(\ref{nonleadvar}).
Also note that from (\ref{leadvar}),  since
$x(G_i)=x_{\ell_i}+d(G^S_i),$
we have
$$x(G_i)=f(G^{(i)})-f(G^{(i-1)}), $$
which, along with (\ref{g1a}), leads to the relations in (\ref{g2}).

Next, we show $x(A)\le f(A)$, for any $A\subseteq E$.
Write $A= \cup_{i=1}^m D_i $, with $D_i \subseteq G_i$, $i=1,...,m$.
Similar to the $G^{(i)}$ notation, denote
$D^{(i)}:=\cup_{k=1}^i D_k$, and $D^{(0)}:=\emptyset$.
Then, from (\ref{g1a}), we have
\begin{eqnarray*}
x(D_i) &\leq & f(G^{(i-1)},D_i) - f(G^{(i-1)}) \\
&\le& f(D^{(i)}) - f(D^{(i-1)})
\end{eqnarray*}
where the second inequality follows from submodularity. Hence,
summing over all $i$, taking into account that $D_i$'s are non-overlapping
(since the $G_i$'s are), we have
\begin{eqnarray*}
x(A) =\sum_{i=1}^m x(D_i) \le f(D^{(m)})= f(A).
\end{eqnarray*}
\hfill$\Box$

\begin{thm}
\label{thm:optimal}
{\rm
 The feasible solution characterized in Lemma \ref{lem:feasible}
is optimal if the indices of the lead and nonlead variables
satisfy the following relations:
\begin{eqnarray}
\label{optc}
{\ell_j}\le {\ell_{j+1}}, \; j=1,...,m-1; \qquad
{\ell_j} \le {k_j}, \quad k_j\in G^S_j, \;j=1,...,m.
\end{eqnarray}
}
\end{thm}

{\startb Proof.}
The dual of the LP in (\ref{lpside}) takes the following form:
\begin{eqnarray}
\label{dual}
\min && \sum_{A\subseteq E} y_A f(A) -\sum_{i\in E} d_iz_i \\
{\rm s.t.} &&\sum_{A\ni i} y_A -z_i \ge c_i, \qquad i\in E; \nonumber\\
&& y_A\ge 0, \quad z_i\ge 0; \qquad A\subseteq E, \; i\in E; \nonumber
\end{eqnarray}
where $y_A$ and $z_i$ are the dual variables (the latter corresponds
to the side constraints).

Following (\ref{nonleadvar}) and (\ref{leadvar}),
the primal objective value is:
\begin{eqnarray}
\label{primalobj}
\sum_{j=1}^m (c_{\ell_j}-c_{\ell_{j+1}}) f(G^{(j)}) -
\sum_{j=1}^m \sum_{k_j\in G^S_j} (c_{\ell_j}-c_{k_j})d_{k_j},
\end{eqnarray}
where $c_{\ell_{m+1}}:=0$.

Now, let the dual variables take the following values:
$$y_A=c_{\ell_j}-c_{\ell_{j+1}}, \qquad A=G^{(j)}, \quad j=1,...,m;$$
let all other $y_A=0$; let
$$z_{k_j}=c_{\ell_j}-c_{k_j}, \qquad k_j\in G^S_j, \quad j=1,...,m;$$
and let all other $z_i=0$. Then, the non-negativity of these dual
variables follows from (\ref{optc}), taking into account
(\ref{corder}). Furthermore, consider the $h$-th constraint in the
dual problem. Suppose  $h\in G_j$. If $h$ is a lead variable, i.e.,
$h=\ell_j$, then $h\in G^{(i)}$ for all $i\ge j$. Hence the
constraint becomes:
$$\sum_{i\ge j} (c_{\ell_i}-c_{\ell_{i+1}}) = c_{\ell_j} = c_h.$$
If $h$ is a non-lead variable, i.e., $h=k_j\in G^S_j$, then we have
$$z_h=z_{k_j}=c_{\ell_j}-c_{k_j}=c_{\ell_j}-c_h,$$
as well as $h\in G^{(i)}$ for all $i\ge j$.
In this case,  the constraint becomes:
$$\sum_{i\ge j} (c_{\ell_i}-c_{\ell_{i+1}}) - (c_{\ell_j}-c_h) = c_h.$$

The complementary slackness condition is also readily checked.
(In particular, the dual objective is equal to the primal objective
in (\ref{primalobj}).)

Finally, primal feasibility has already been established in
Lemma \ref{lem:feasible}.
Hence, the optimality of both primal and dual solutions.
\hfill$\Box$

\section{The Grouping Algorithm}
\label{sec:algorithm1}

The grouping algorithm to be presented below forms the groups to
maintain feasibility. In particular it maintains the distinction of
non-lead variables and lead variables in (\ref{nonleadvar}) and
(\ref{leadvar}), and enforces the two conditions in Lemma
\ref{lem:feasible} --- that the lead variable of each group
satisfies the side constraint, and that (\ref{g1a}) holds for each
group. Furthermore, we introduce a grouping structure\footnote{Lack
a better choice, we name them subgroups. However, there is no
connections between our problem and the algebra entities of group
and subgroup} within each group to form a hierarchical structure
that is crucial for the feasibility proof, as well as track the
order of variables entering the group. Optimality follows, as argued
in Lemma \ref{lem:feasible} and Theorem \ref{thm:optimal} above.

Although not essential to understanding the algorithm, it would
help intuition to view the groups as associated with different
priorities, with $G_1$ having the highest priority, and $G_m$ the
lowest priority. The variables can be viewed as throughput
associated with each class, whereas the side constraints represent
the minimal requirements on the throughput. The lead variable in
each group contributes to maximizing the objective function, a
weighted sum of the throughput of all job classes. The non-lead
variables are non-contributors; as such, they are moved into a
higher-priority group parsimoniously just to satisfy their
respective side constraints.

Inductively, assume we know how to do the grouping when $|E|=n-1$.
Specifically, suppose this results in $m-1$ groups,
${G_1,...,G_{m-1}}$. Within each group $G_i$, we have,
\begin{itemize}
\item
lead variable $\ell_i$;
\item
subgroup $g^{i}_j$, $j=1, \cdots, j_i$;
\end{itemize}
and they satisfies, $\forall k_j \in g^i_j$,
\begin{eqnarray}
\label{subgroup_1} f(G^{(i-1)}; g^i_1, \cdots, g^i_{j-1}) + d_{k_j}
\le f(G^{(i-1)}; g^i_1, \cdots, g^i_{j-1}, k_j)
\end{eqnarray}
and there exists at least one $k_\ell \in g^i_j$ such that,
\begin{eqnarray}
\label{subgroup_2} && f(G^{(i-1)}, g^i_1,\cdots, g^i_{j-1}, k_j,
k_\ell) - f(G^{(i-1)},
  g^i_1, \cdots, g^i_{j-1} , k_j) \nonumber \\ &&+ f(G^{(i-1)};g^i_1, \cdots,
  d^i_{j-1}) + g_{k_j} < d_{k_\ell}
\end{eqnarray}
and,
\begin{eqnarray}
\label{leadvariable} f(G^{(i)})-f(G^{(i-1)}) - d(G_i^S) \ge
d_{\ell_i}.
\end{eqnarray}

Now consider the added $n$-th dimension; in particular, we want to
find where to put the new element $n$ and what value to assign to $x_n$.
In the algorithm presented below, this is accomplished in Step 1.
Then, we need to adjust the newly generated group
and the lower priority groups that are affected
(by the new group) so as to maintain
feasibility; and this is done in Steps 2 and 3 of the algorithm.

To start with, we let $n$ itself form the $m$-th group.

\bigskip
\noindent {\bf The Grouping Algorithm}

\medskip
\noindent
{\bf Step 1. ({\it Group Identification})}
\begin{itemize}
\item
Find the largest $i\le m $, such that
\begin{eqnarray}
\label{move}
f(G^{(i-1)}) +d_{n} \le f(G^{(i-1)};  n),
\end{eqnarray}
and
\begin{eqnarray}
\label{positive}
f(G^{(i)})+d_n  >f(G^{(i)},n).
\end{eqnarray}
\end{itemize}

\medskip
\noindent
{\bf Step 2. ({\it Subgroup Identification})}
Within group $G^i$, find the largest $j$, such that
\begin{eqnarray}
\label{subgroup_formation_1} f(G^{(i-1)}, g^i_1, \cdots, g^i_{j-1})
+ d_n  \le f(G^{(i-1)}; g^i_1, \cdots, g^i_{j-1}, n),
\end{eqnarray}
form new subgroups according to the validity of the following
inequality,
\begin{eqnarray}
\label{subgroup_formation_2} f(G^{(i-1)}, g^i_1, \cdots,
g^i_{j-1},n) + d_{k_j}  \le f(G^{(i-1)}; g^i_1, \cdots, g^i_{j-1},
n, k_j),
\end{eqnarray}
\begin{itemize}
\item if for all $k_j\in g^i_j$, (\ref{subgroup_formation_2}) holds,
then, $n$ will becomes new subgroup $g^i_j$, and the indices for all
the previous subgroups that have index larger and equal to $j$ will
have to increase by one;
\item otherwise $n$ will join the subgroup $g^i_j$.
\end{itemize}

\medskip
\noindent {\bf Step 3. ({\it Subgroup Validation})} This step is
to determine if there is a feasible solution that the variables in
subgroup $g^i_j$ are valued at the side constraints. This is
equivalent to solve the membership problem for vector $(d_{k_1},
d_{k_2}, \cdots, d_{k_\ell})$ to the polymatroid defined by the
following submodular function,
\begin{eqnarray}
\label{sub_polymatroid} g(A) := f(G^{(i-1)}, g^i_1,\cdots,
g^i_{j-1}, A) - f(G^{(i-1)};g^i_1, \cdots,
  d^i_{j-1})
\end{eqnarray}
i.e. the subgroup is valid if the following LP has feasible
solution,
\begin{eqnarray}
\label{mp_LP} \min x_{k_1}, s.t. {\bf x} \in P(g), x_{k_i} \ge
d_{k_i}
\end{eqnarray}
If the LP in (\ref{mp_LP}) is founded to be infeasible, then this
subgroup is invalid, hence to be dissolved. Then if $j>1$, then we
need to go back to Step 2 to identify the subgroup for all the
member previously in $g^i_j$; if $j=1$, the group $G^i$ need to be
dissolved, go back to Step 1 for group identification for each
member previously in $G^i$.
\medskip
\noindent {\bf Step 4. }
\begin{itemize}
\item
Conduct validation step for subgroup $g^i_{j+1}, g^i_{j+2}$;
\item
Repeat the procedure to all the variables that has been bypassed;
\end{itemize}

\begin{thm}
\label{thmpb} {\rm The grouping algorithm generates a solution that
follows (\ref{nonleadvar}) and (\ref{leadvar}), and satisfies both
feasibility requirements (a) and (b) of Lemma \ref{lem:feasible},
and the inequalities in (\ref{optc}) of Theorem \ref{thm:optimal}.
Hence, it is an optimal solution to the LP in (\ref{lpside}). }
\end{thm}

\noindent {\bf Proof.} We only need to check the feasibility. To do
this, we need to check if our output satisfies equation (\ref{g1a}).
For any $A\subset G_i$, let us consider two cases, (i) $\ell_i \in
G_i-A$; (ii) $\ell_i\in A$.

\noindent {\bf Case I:} Certainly, we have, $A= \cup_{j=1}^{j_i}
A_j$, with $A_j= A\cap g^i_j$. From the validation step, we know
that for each $j$,
\begin{eqnarray*}
x(A_j) =d(A_j) &\le&  f(G^{(i-1)}, g^i_1,\cdots,g^i_{j-1},
A_j)-f(G^{(i-1)}; g^i_1,\cdots,g^i_{j-1})
\end{eqnarray*}
Due to submodularity, we have,
\begin{eqnarray*}
x(A_j) =d(A_j) &\le&  f(G^{(i-1)},  A_1,\cdots,A_{j-1},
A_j)-f(G^{(i-1)}; A_1,\cdots,A_{j-1})
\end{eqnarray*}
It is true for for each $j$, then, sum them up, it implies that
(\ref{g1a}) is true.

\noindent {\bf Case II:}


For any subgroup $g_j^i$, let us denote $k_1^j, k_2^j,\cdots,k_h^j$
of its member according to the order they entered the subgroup.
Hence, the inequality (\ref{positive}) must be satisfied in the
following form,
\begin{eqnarray*}
d_{k^j_1} &>&
f(G^{(i-1)};g_1^i, \cdots, g^i_{j-1}; \ell_i, k^j_1) -f(G^{(i-1)} ; g_1^i, \cdots, g^i_{j-1};\ell_i), \\
k_{g^j_2} &>&
f(G^{(i-1)}; g_1^i, \cdots, g^i_{j-1};\ell_i, k^j_1, g^i_2) -f(G^{(i-1)}; g_1^i, \cdots, g^i_{j-1};\ell_i, k^j_1), \\
&\cdots\cdots& \\
d_{k^j_h} &>& f(G^{(i-1)} ; g_1^i, \cdots, g^i_{j-1}; \ell_i, k^j_1,
\dots k^j_{h-1}, k^j_h) -f(G^{(i-1)} ; \ell_i, g^i_1, \dots
g^i_{h-1}) ;
\end{eqnarray*}
Because when $k_1^j$ enters the group, $G{(i)}$ is made of
$G^{(i-1)}$, subgroups $g_1^i, \cdots, g^i_{j-1}$ and its lead
member $\ell_i$, so the inequality (\ref{positive}) takes the form
of the first inequality, and so on.

Meanwhile,
\begin{eqnarray*}
&&
f(G^{(i-1)};g_1^i, \cdots, g^i_{j-1}; \ell_i, k^j_1) -f(G^{(i-1)} ; g_1^i, \cdots, g^i_{j-1};\ell_i), \\
&&
f(G^{(i-1)}; g_1^i, \cdots, g^i_{j-1};\ell_i, k^j_1, k^i_2) -f(G^{(i-1)}; g_1^i, \cdots, g^i_{j-1};\ell_i, k^j_1), \\
&&\cdots\cdots \\
&& f(G^{(i-1)} ; g_1^i, \cdots, g^i_{j-1}; \ell_i, k^j_1, \dots
k^j_{h-1}, k^j_h) -f(G^{(i-1)}; \ell_i, g^i_1, \dots g^i_{h-1}) ;
\end{eqnarray*}
forms a solution to another polymatroid, from its property, we can
conclude that,
\begin{eqnarray*}
d(A_j) > f(G^{(i-1)};g_1^i, \cdots, g^i_{j-1}, g^i_j; \ell_i)
-f(G^{(i-1)} ; g_1^i, \cdots, g^i_{j-1};A_j^c, \ell_i)
\end{eqnarray*}
where $A^c:=G^S_i-A$ and $A^c_j = A^c\cap g^i_j$.



Therefore, for any $A\subseteq G^S_{i+1}$, we have
\begin{eqnarray}
\label{more} d(A) >  f(G^{(i)}; A, \ell_{i+1})- f(G^{(i)};
\ell_{i+1}).
\end{eqnarray}
\begin{eqnarray*}
&&x_{\ell_{i+1}}+ d(A)\\
&=& f(G^{(i)};G_{i+1})- f(G^{(i)}) - d(A^c)\\
&\le & f(G^{(i)};G_{i+1})- f(G^{(i)}) -
 [f(G^{(i)}; A^c, \ell_{i+1})- f(G^{(i)}; \ell_{i+1})]\\
&\le &  f(G^{(i)}; A, \ell_{i+1})- f(G^{(i)}),
\end{eqnarray*}
where the first inequality follows from (\ref{more}), applied to
$A^c$, and the second one follows from submodularity, i.e.,
\begin{eqnarray*}
f(G^{(i)};G_{i+1}) -f(G^{(i)}; A^c, \ell_{i+1})
\le   f(G^{(i)}; A, \ell_{i+1})- f(G^{(i)}; \ell_{i+1});
\end{eqnarray*}
hence, the feasibility.

\hfill$\Box$

Consider the $b$-version of the performance polytope (in Definition
\ref{def:cl}) with side constraints:
\begin{eqnarray}
\label{bside}
\calp'(b) = \{
x(A) \ge b(A), A\subseteq E;\,
 x_i \le d_i , i\in E\}.
\end{eqnarray}
(For instance, $d_i$ is the upper limit on the delay of type $i$
jobs.) We want to solve the following LP problem:
\begin{eqnarray}
\label{lpside1} \min_{x\in\calp'(b)}  \sum_{i=1}^n c_i x_i.
\end{eqnarray}
Direct verification will establish that the Grouping Algorithm is
readily adapted to solve the above problem: simply change $f$ to
$b$, and reverse the directions of all the inequalities involved.
(Note that the side constraints are now $x_i\le d_i$ instead of
$x_i\ge d_i$ as in the $f$-version.)

\subsection{The Complexity of the Grouping Algorithm}

In general, the membership problem in {\bf Step 3} is as hard as the
submodular function minimization itself. However, the special
structure obtained through our construction of the subgroup will
enable us to solve this problem much easily. The key is the
following important observation,
\begin{lem}
In (\ref{mp_LP}), only the side constraints and the constraint for
the full set can be tight.
\end{lem}
{\bf Proof} We prove this by induction. The case of two variables is
obviously trivial. To illustrate the main idea, we first consider
the case of three, then the general case. Now, suppose that subgroup
$g^i_j$ consists of two elements, $k_1$ and $k_2$. According to the
{\bf Step 2} of the grouping algorithm, we have,
\begin{eqnarray}
\label{eqn:2_var_1}
 f(G^{(i-1)};g^i_1, \cdots, d^i_{j-1}) + d_{k_j} \le
f(G^{(i-1)}; g^i_1, \cdots, g^i_{j-1}, k_j), j=1,2
\end{eqnarray}
\begin{eqnarray}
\label{eqn:2_var_2} && f(G^{(i-1)}, g^i_1,\cdots, g^i_{j-1}, k_j,
k_\ell) - f(G^{(i-1)},
  g^i_1, \cdots, g^i_{j-1} , k_j) < d_{k_\ell}, j, \ell =1,2.
\end{eqnarray}
Then, adding variable $k_3$ that satisfies,
\begin{eqnarray*}
f(G^{(i-1)};g^i_1, \cdots, g^i_{j-1}) + d_{k_3}  \le f(G^{(i-1)};
g^i_1, \cdots, g^i_{j-1}, k_3),
\end{eqnarray*}
and
\begin{eqnarray*}
f(G^{(i-1)};g^i_1, \cdots, g^i_{j}) + d_{k_3}  > f(G^{(i-1)}; g^i_1,
\cdots, g^i_{j}, k_3).
\end{eqnarray*}
First, we know that,
\begin{eqnarray*}
 f(G^{(i-1)};g^i_1, \cdots, d^i_{j-1}, k_3) + d_{k_j} \le
f(G^{(i-1)}; g^i_1, \cdots, g^i_{j-1}, k_3, k_j), j=1,2
\end{eqnarray*}
can not be be valid for both $j=1,2$, since in this case, $k_3$ will
form a subgroup of its own according to Step 2 of the grouping
algorithm. From the fact that we have at least one $j^*$, such that,
\begin{eqnarray*}
 f(G^{(i-1)};g^i_1, \cdots, d^i_{j-1}, k_3) + d_{k_{j^*}} <
f(G^{(i-1)}; g^i_1, \cdots, g^i_{j-1}, k_3, k_{j^*}),
\end{eqnarray*}
it implies that $x_{k_3}\le g(k_3) $ can not be tight, for otherwise
the solution is infeasible. Now consider the subsets of $\{k_1,
k_3\}$ and $\{k_2, k_3\}$, they can not be tight, either. Otherwise,
suppose that $\{k_1, k_3\}$ is tight, then
\begin{eqnarray*}
x_{k_2}\le g(k_1, k_2, k_3)- g(k_1, k_3) \le g(k_1, k_2)- g(k_1) <
d_2
\end{eqnarray*}
which leads to infeasibility of the problem. Here the first
inequality is due to submodularity and the second one is just
(\ref{eqn:2_var_2}).  Similar arguments apply to $\{k_2, k_3\}$.

For the general $n$ case, we can repeat the similar argument.
Suppose that we have a $k_n$ added into the subgroup $g_i^j$. Then
again, due to the fact that $k-n$ is admitted into the subgroup, the
following inequality will hold for some existing member of $g^i_j$,
$j^*$,
\begin{eqnarray*}
 f(G^{(i-1)};g^i_1, \cdots, d^i_{j-1}, k_n) + d_{k_{j^*}} <
f(G^{(i-1)}; g^i_1, \cdots, g^i_{j-1}, k_n, k_{j^*}),
\end{eqnarray*}
which implies that $x_{k_n}\le g(k_n) $ can not be tight. For all
the subset contains $k_n$, but not the full set,  the above
submodularity argument can be employed again to show that the
corresponding inequality can not tight. For those do not contain
$k_n$, the induction hypothesis tells us that they can not be tight.
Hence, the proposition follows. $\Box$

\begin{lem}
The complexity of {\bf Step 3} is $O(LP(n))$. Where $LP(n)$ denote
the complexity of solving a linear program with no more than $n$
constraints and $n$ variables.
\end{lem}
{\bf Proof } By complementary slackness, in the dual of
(\ref{mp_LP}), that is,
\begin{eqnarray}
\label{mp_LP_D} \max \sum_A y^Ag(A)-\sum_{k_i} d_{k_i}z_{k_i},
s.t. \sum_{A \ni k_1} y^A -z_{k_1}\le 1, \sum_{A \ni k_i} y^A
-z_{k_i}\le 0, \forall k_i\neq k_1,
\end{eqnarray}
the variables that corresponds to constraints other than the side
constraints and the one for the full set equal to zero. Therefore,
to solve the dual problem, we only need to solve a linear program
with $\ell+1$ variables and $\ell$ constraints. That gives us the
complexity bound of $O(LP(n))$. In fact, to check whether the
original problem is feasible, we simply need to rule out the
unboundness of the dual problem. $\Box$

\medskip

To characterize the computational effort of the entire grouping
algorithm, consider the point where we started Step 1 of the
algorithm, i.e., when the $n$-th element is added to the $m-1$
groups that have already been formed for the problem of $n-1$
elements. Note that throughout the algorithm, whenever an element
moves, individually or in group, it always moves into a
higher-priority (lower-indexed) group, or to a higher-priority
within a group; and the very first subgroup is invalid, an entire
group is eliminated. (Indeed, the number of groups can only
decrease.) Hence, for an element originally in the $k$-th group, the
number of moves can be no greater than total number of subgroups
ahead, with each move involving a fixed number of operations
associated with checking the inequalities in
(\ref{move}),(\ref{positive}),({\ref{subgroup_formation_1}), and
({\ref{subgroup_formation_2})), and solving an LP. Therefore, to
finalize the grouping with the additional $n$-th element, the
computational effort involved is $O(n^2\times LP(n))$. Since we
(inductively) bring in the $n$ elements one at the time, the overall
effort of the grouping algorithm is $O(n^3\times LP(n))$.


\medskip

\subsection{Connections to the Queueing Scheduling Problem}

Return to the queueing scheduling problem considered in
\S\ref{sec:qnet}), with the side constraints.

\begin{thm}
\label{thm:qschedside}
{\rm
Consider the control problem
in (\ref{qsched}), with additional side constraints,
$x_i\ge d_i$, $i\in E$, and suppose the $c_i$'s follow (\ref{corder}).
Solve the LP in (\ref{lpside}).
Let $G_1,\dots, G_m$ be the groups formulated by applying
the grouping algorithm.
Then,  the optimal control $u^*$ is a priority rule
among the groups, with $G_1$ having the highest priority,
and $G_m$ the lowest priority; and within each group,
$u^*$ is a randomized rule that enforces the side constraints
(on the non-lead variables).
}
\end{thm}

For the last statement in the above theorem, reason as follows.
With side constraints, the optimal solution $x^*$
to the LP in (\ref{lpside}) returned by
the grouping algorithm is an interior point of the original
polymatroid $\calp (f)$ (i.e., without the side constraints).
Following Caratheodory's theorem (refer to, e.g.,
Chv\'atal \cite{chvatal}), $x^*$
can be expressed as a convex
combination of no more than $n$ vertices
of $\calp (f)$, each of which corresponds
to a priority rule.
(Note that the equality $x(E)=f(E)$ has reduced the
degree of freedom by one; hence, $n$, instead of $n+1$.)
In other words, $x^*$ can be
realized by a control that is a {\it randomization} of at most
$n$ priority rules, with the coefficients of the convex combination
used as the probabilities for the randomization.


In terms of implementing the optimal control,
one can do better than randomization. It is known
(e.g., Federgruen and Groenevelt \cite{fgmgc})
that any interior point of the
performance space can be realized by a particular
dynamic scheduling policy, due originally to
Kleinrock \cite{kleinrock,kleinrock2},
in which the priority index of each job present in the system
grows proportionately to the time it has spent waiting in queue,
and the server always serves the job that has the highest index.
This scheduling policy is completely specified by the
proportionate coefficients associated with the jobs classes,
which, in turn, are easily determined by the performance vector
(provided it is at the interior of the performance space).
In terms of practical implementation, there are several versions
of this scheduling policy, refer to
\cite{timefunction,delaycost}.

\section{Submodular Function Minimization}
\label{sec:minsubmod}

Suppose we want to minimize a submodular function
$\psi(A)$, $A\subseteq E$.
Following  \cite{cunningham}, we can transform
this into an equivalent problem as follows.
Define
$$d_i := \psi(E-\{i\}) -\psi(E), \qquad i\in E;$$
and write $d(A)=\sum_{i\in A}d_i$. Let
\begin{eqnarray}
\label{fa}
 f(A) :=\psi(A) + d(A) -\psi(\emptyset).
\end{eqnarray}
Then, minimizing
$$f(A) - d (A)=\psi(A)  -\psi(\emptyset)$$
is equivalent to minimizing the original $\psi(A)$.
Clearly, $f(\emptyset) =0$, and $f$ is submodular.
To see that $f$ is also increasing, consider $A\subset E$,
and $i\not\in A$. We have:
\begin{eqnarray*}
f(A+\{i\})-f(A)
&=&
\psi(A+\{i\})+d(A)+d_i-\psi(A)-d(A)\\
&=&
\psi(A+\{i\})-\psi(A)
+\psi(E-\{i\})-\psi(E)\\
&\ge& 0,
\end{eqnarray*}
where the inequality follows from the submodularity of $\psi$.

Hence, the minimization of a submodular function $\psi(A)$
over $A\subseteq E$ is equivalent to the following problem:
\begin{eqnarray}
\label{minf}
\min_{A\subseteq E} \{ f(A) - d(A)\}.
\end{eqnarray}
To solve this problem using the grouping algorithm
of the last section, we still need (i) $d_j\ge 0$ for all $j\in E$,
and (ii) $f(A)\ge d(A)$ for all $A\subseteq E$
(refer to Lemma \ref{lemrossyao}).

Consider (i) first. If $d_j<0$, for some $j$,
then, for any set $A\subset E$, due to the
submodularity of $\psi$, we have, for any
negative $d_j$,
\begin{eqnarray}
\label{ignore}
0>d_j=\psi (E- \{j\})-\psi (E) \ge
\psi (A- \{j\})-\psi (A),
\end{eqnarray}
implying $\psi (A- \{j\})>\psi (A)$.
That is, any minimizing set $A$ will not contain $j$ if $d_j<0$.
Hence, we can redefine the set $E$ by first excluding
such $j$'s. This way, we can effectively assume $d_j\ge 0$
for all $j$.


The second requirement,
$f(A)\ge d(A)$,
following (\ref{fa}),
is equivalent to  $\psi (A)\ge \psi (\emptyset )$ for all $A$.
To overcome this difficulty,
we shall assume $\psi (A)$ is bounded from below:
there exists an $M>0$ such that $\psi (A)\ge -M$
for all $A\subseteq E$.
Then, we can replace $\psi$ by $\psi'$ defined as follows:
$$\psi'(A)=\psi (A), \quad  A\neq \emptyset ; \qquad
\psi'(\emptyset ) =-M;$$
and $\psi'(A)\ge\psi'(\emptyset)$ for all $A\subseteq E$.
We then solve the minimization problem in (\ref{minf}),
but over $E-\emptyset$ (instead of $E$).
Once the minimizing set $A^*$ is identified, 
we can compare $\psi(A^*)$ with $\psi (\emptyset)$ to determine whether
$A^*$ or $\emptyset$ is optimal.

\begin{thm}
\label{inclexcl}
{\rm
Let $G_1$ be the first group formulated from
applying the grouping algorithm (under the permutation
$\pi =(1,...,n)$).
(Recall, $x_1$ is the lead variable in $G_1$.)
Then, for any
 $B\neq G_1$, with $1\in B\subseteq E$, $B$ is inferior to
$G_1$ in terms of
yielding a larger objective value in (\ref{minf}).
%
}
\end{thm}

{\startb Proof.}
First, suppose $B\subset G_1$.
From Lemma \ref{lem:feasible}, we have
$$f(B)\ge x(B), \qquad f(G_1)=x(A),$$
where $x$ denotes the variable of the LP over $\calp' (f)$
as in the last section.
Hence,
$$f(G_1)- f(B)\le x(G_1)-x(B)
=d(G_1)-d(B),$$
where the equality follows from the fact that
$1\in B\subset G_1$;
hence, $G_1-B$ only contains non-lead variables.

Next, suppose $B\supset G_1$.
Write  $B=G_1+A$, where
 $A\subseteq E-G_1$.
We have, following Lemma \ref{lempb},
$$f(B)=f(A+G_1) \ge x(A)+x(G_1) =x(A)+f(G_1).$$
Hence,
$$f(B)-f(G_1)\ge x(A)\ge d(A)=d(B)-d(G_1).$$
%

Now, consider a general $B$. From submodularity, we have
$$f(G_1)+f(B) \ge f(G_1\cup B)+f(G_1\cap B).$$
Hence,
\begin{eqnarray*}
&&f(B)-d(B)\\
&\ge& f(G_1\cup B)+f(G_1\cap B)-f(G_1)-d(B)\\
&=& f(G_1\cup B)- d(G_1\cup B)
+f(G_1\cap B)-d(G_1\cap B)-f(G_1)+d(G_1)\\
&\ge& 2[f(G_1)-d(G_1)]-f(G_1)+d(G_1)\\
&=&f(G_1)-d(G_1) ,
\end{eqnarray*}
where the first equality follows from $d(\cdot)$ being
additive, and the second inequality follows from
the two cases already established above (regarding the subsets
and supersets of $G_1$).
\hfill$\Box$

Following Theorem \ref{inclexcl}, we can conclude that
after the set $A=G_1$ is identified
we do not have to consider the element $1$ any more: any set that
includes it must yield a larger objective value in (\ref{minf})
than $G_1$. To obtain the solution to problem(\ref{minf}), let us
add another variable $\{0\}$ into the problem, define $d_0=-1$ and
$f(A\cup\{0\}) =f(A)$, this again defined a submodular function on
${\hat E}=E\cup \{0\}$. In the running of the algorithm, give $0$
the highest priority, then it becomes the lead variable in $G_1$,
hence the algorithm will produce,
$$\min_{A\subseteq {\hat E}-\{0\}} \{ f(A) -d(A)\}.$$
However, from the definition, we know that $0$ will never belongs
to the minimum, hence, this is the optimum of (\ref{minf}).

\begin{rem}
\label{rem:ours}
{\rm
Another way to view Theorem \ref{inclexcl} is to consider the
minimization problem in (\ref{minf}), with the restriction
that $A\ni 1$, i.e.,
\begin{eqnarray}
\label{minf1}
\min_{1\in A\subseteq E} \{ f(A) -d(A-\{1\})\}.
\end{eqnarray}
Note that we have dropped the $-d_1$ term in (\ref{minf1}).
The problem in (\ref{minf1})
can be related to the LP over $\calp' (f)$, with the following
cost coefficients:
$$c_1=1, \quad c_2=\cdots=c_n=0.$$
Then, (\ref{minf1}) becomes the same as the dual LP
problem in (\ref{dual}), with $y_{G_1}=1$, and
$y_A =0$ for all other $A=G^{(j)}$, $j=2,...,m$;
$z_{k_j}=1$, $k_j\in G^S_1$, and all other $z_i=0$.
Hence, the grouping algorithm will return
$A=G_1$ as the optimal solution to (\ref{minf1}).

This also suggests an alternative approach to the SFM in (\ref{minf}):
for $i=1,\dots ,n$, solve $\min_{i\in A\subseteq E} \{ f(A) -d(A-\{i\})\}$,
 by applying the grouping algorithm to the
LP over $\calp' (f)$, with $c_i=1$ and $c_j=0$ for $j\neq i$.
Each results in a subset $G_1\ni i$. Clearly, only these $n$ subsets
need to be considered; hence, we can compare these in terms of
their corresponding objective values in (\ref{minf}), and pick the
the smallest one.
\hfill$\Box$
}
\end{rem}

\begin{rem}
\label{rem:others}
{\rm
The starting point of virtually all existing SFM algorithms
in the literature
(in particular, those overviewed in the introductory section)
is to relate the problem in (\ref{minf}) to the following
equivalent LP:
\begin{eqnarray}
\label{otherLP}
\max \{x(E): \; x\in\calp (f),\; x_i\le d_i, \, i\in E\},
\end{eqnarray}
where $\calp (f)$ is the polymatroid in (\ref{defPf}).
(The equivalence between (\ref{minf}) and (\ref{otherLP}) is
established in \cite{edmonds}; also refer to \cite{mccormick}.)
The above LP is then solved via finding augmenting paths in
an associated graph. In this approach, the key step
is to make sure that any improvement on the current solution
$y$ remains within the (base) polymatroid $\calb (f)$, and this is accomplished via
Carath\'eodory's representation of $x$ (as a convex combination
of the vertices of the polymatroid).
As commented by McCormick in \cite{mccormick}, ``this is a rather
brute-force way to verify that $y\in\calb (f)$.
However, 30 years of research have not yet
produced a better idea.''

In contrast, as explained in Remark \ref{rem:ours},
our approach is to solve the SFM in (\ref{minf}) via
$n$ subproblems such as the one in (\ref{minf1}), which corresponds
to the dual of an LP over $\calp' (f)$. Note that $\calp' (f)$
is different from the polytope in (\ref{otherLP}) --- the inequalities
in the side constraints are in the opposite direction.
This subtle distinction
turns out to be crucial: While the grouping algorithm
solves the LP over $\calp' (f)$, it is not clear how the algorithm
can be adapted to solve the LP in (\ref{otherLP}).
Recall, the intuition behind the grouping algorithm is that each
job class associated with a non-lead variable,
a non-contributor to maximizing the objective
function, is moved into a
higher-priority group just to meet the {\it lower} limit
on its performance (throughput) requirement. Should the side
constraints become upper limits, the structure of the solution would have to
be completely different.
\hfill$\Box$
}
\end{rem}

\subsection*{Acknowledgements}

We are much indebted to Shuzhong Zhang for
his detailed and insightful comments
on an earlier version of the paper.
We also thank Jay Sethuraman and Li Zhang for
useful discussions at
various stages of this research.

\end{document}